\documentclass[times,review,sort&compress,11pt]{elsarticle}
\usepackage{epstopdf}
\usepackage[caption=false]{subfig}
\usepackage[hidelinks]{hyperref}
\usepackage[active]{srcltx}
\usepackage{amsthm}
\usepackage{amssymb}
\usepackage{stackrel}
\usepackage{xcolor}
\usepackage{amsmath,amscd}

\newtheorem{dl}{Theorem}[section]
\newtheorem{tl}[dl]{Corollary}
\newtheorem{yl}[dl]{Lemma}

\newtheorem{lz}[dl]{Example}
\newtheorem{wenti}[dl]{Nahm's conjecture}

\newtheorem{remark}[dl]{Remark}

\numberwithin{equation}{section}

\newcommand{\be}{\begin{equation}}
\newcommand{\ee}{\end{equation}}
\newcommand{\ba}{\begin{array}}
\newcommand{\ea}{\end{array}}
\newcommand{\bmn}{\begin{eqnarray}}
\newcommand{\emn}{\end{eqnarray}}
\newcommand{\bnm}{\begin{eqnarray*}}
\newcommand{\enm}{\end{eqnarray*}}
\newcommand{\bln}{\begin{subequations}}
\newcommand{\eln}{\end{subequations}}

\newproof{pot333}{\bf Proof of Theorem \ref{wanngliuquan}}
\newproof{pot444}{\bf Proof of Theorem \ref{mainthm}}
\newproof{pot777}{Proof of Theorem \ref{secondeq-third}}

\newproof{pot555}{Applications of Theorem \ref{firsteq}}
\newproof{pot666}{Applications of Theorem \ref{secondeq-second}}
\newproof{pot888}{Applications of Theorem \ref{secondeq-third}}

\def\qed{\hfill \rule{4pt}{7pt}}
\def\pf{\noindent {\it Proof.} }

\begin{document}

\title{A general $q$-series transformation and its applications to multi-sum   Rogers-Ramanujan-Slater identities}
\author{Jianan Xu\fnref{fn1}}
\fntext[fn1]{E-mail address: 20224007007@stu.suda.edu.cn}
\address[P.R.China]{Department of Mathematics, Soochow University, SuZhou 215006, P.R.China}
\author{Xinrong Ma\fnref{fn2,fn3}}
\fntext[fn2]{This  work was supported by NSFC grant No. 11971341}
\fntext[fn3]{Corresponding author. E-mail address: xrma@suda.edu.cn.}
\address[P.R.China]{Department of Mathematics, Soochow University, SuZhou 215006, P.R.China}
\markboth{J.N. Xu and  X. R. Ma}{A general $q$-series transformation and its applications}
\begin{abstract}In the present  paper, we establish a general transformation for $q$-series which contains L. Wang et al's transformation involved in Nahm series. As direct applications, some concrete new transformation formulas for the ${}_{r+1}\phi_r$ series as well as $q$-identities of multi-sum   Rogers-Ramanujan-Slater type are presented.\end{abstract}
\maketitle
\vspace{20pt}
\parskip 7pt

\section{Introduction}
It is well known that  all sum-product $q$-series identities as the form
\begin{align}
\sum_{n=0}^{\infty} \frac{q^{n^2}}{(q ; q)_n}=\frac{1}{\left(q, q^4 ; q^5\right)_{\infty}}, \quad \sum_{n=0}^{\infty} \frac{q^{n(n+1)}}{(q ; q)_n}=\frac{1}{\left(q^2, q^3 ; q^5\right)_{\infty}}\label{rogres}
\end{align}
are called to be   Rogers-Ramanujan type identities since it is  Rogers and Ramanujan who discovered \eqref{rogres} independently around 1900. As of today,    such Rogers-Ramanujan type identities serve as one of the witnesses for deep connections between the theory of $q$-series and modular forms. Nevertheless, the problem of  determining  certain given $q$-hypergeometric series is modular or not is
    interesting but not easy to solve. Motivated by applications in conformal field theory, W. Nahm \cite{name} paved a new way  to study this question via a special class of multi-sum ple series. His idea can be summarized as follows.
   \begin{wenti}[cf. \mbox{\cite[Sect. 3]{name}}]
   Given a positive integer $r$,  to determine all $r \times r$ positive definite rational matrices $A$, rational vectors $B$ of length $r$ and rational scalars $C$ such that
\begin{align}
f_{A, B, C}(q):=\sum_{n=\left(n_1, \ldots, n_r\right)^{\mathrm{T}} \in(\mathbb{Z}_{\geq 0})^r} \frac{q^{\frac{1}{2} n^{\mathrm{T}} A n+n^{\mathrm{T}} B+C}}{(q ; q)_{n_1} \cdots(q ; q)_{n_r}}
\end{align}
is a modular function. If true, the $(A,B,C)$ is called ``modular triples" and the series
$f_{A, B, C}(q)$ is called ``Nahm sum" and $r$ is ``rank". Hereafter,  $\mathbb{Z}$ denotes the set of all integers and the superscript $T$ stands for the matrix transpose.
\end{wenti}

Up to now, for the rank $r=1$, D. Zagier \cite{13} proved that there are exactly seven modular triples, including two that correspond to the Rogers-Ramanujan identities \eqref{rogres}. When the rank $r=2$, he made an extensive computer search that led to eleven sets of possible modular triples recorded as Table 2 in \cite{13}. All  of these examples except one have been confirmed. Regarding this topic, we refer the reader to L. Wang's paper \cite{12-2} for a comprehensive discussion of all the eleven examples, wherein a full bibliography  can be found. As L. Wang observed, among these examples the fifth and the tenth are quite special. In particular,  all of the eleven examples have been proved except for the fifth example, also there  is no  $q$-series proof for the tenth example. To solve this open problem, in their latest paper \cite{12}, Z. Cao, H. Rosengren, and L. Wang  established the following  $q$-series transformation
 \begin{dl}[\mbox{cf. \cite[Thm. 1.1]{12}}]\label{wanngliuquan}
Let $\left\{A_n\right\}_{n\in \mathbb{Z}}$ be a general sequence such that the series below converge absolutely. Then it holds
\begin{align}
\sum_{i, j \geq 0} \frac{A_{i-j} q^{\frac{j(j-1)}{2}} x^j}{(q ; q)_i(q ; q)_j}=(-x ; q)_{\infty} \sum_{i, j \geq 0} \frac{A_{i-j} q^{\frac{j(j-1)}{2}+i j} x^j}{(q ; q)_i(q ; q)_j(-x ; q)_j}.\label{WWWWWW}
\end{align}
\end{dl}
With the help of \eqref{WWWWWW}, they succeeded to find a $q$-series proof for the fifth example, while this proof can be regarded as the first $q$-series proof for the tenth example. As a crucial step, L. Wang et al's transformation above depends on a special case of the $q$-Chu-Vandermonde summation formula, written in \cite[p. 20]{9-0} by G. E. Andrews as
$$
\frac{1}{(q ; q)_i(q ; q)_j}=\sum_{k=0}^{\min (i, j)} \frac{q^{(i-k)(j-k)}}{(q ; q)_k(q ; q)_{i-k}(q ; q)_{j-k}} .
$$

A quick glance over Transformation \eqref{WWWWWW} inspires us to ask whether there is other similar result in the field of basic hypergeometric series. As far as we are
aware, there are indeed a few  general  transformations found in \cite{10} but no any result like  \eqref{WWWWWW}. For purpose of comparison, we record two of them here as below:

\begin{description}
\item[(i) (\mbox{cf. \cite[(2.2.2)]{10}})]
\begin{align*}
\sum_{k=0}^n \frac{\left(b, c, q^{-n} ; q\right)_k}{(q, a q / b, a q / c ; q)_k} A_k&=\sum_{n\geq i+j,i,j\geq 0} \frac{\left(a q / b c, a q^j, q^{-n} ; q\right)_j}{(q, a q / b, a q / c ; q)_j}\frac{q^{-i j}}{\tau(j)}\frac{\left(q^{j-n}, a q^{2 j} ; q\right)_i}{\left(q, a q^j ; q\right)_i} \left(\frac{b c}{a q}\right)^{i+j} A_{i+j}.
\end{align*}
\item[(ii) (\mbox{cf. \cite[(2.8.2)]{10}})]
\begin{align*}
&\sum_{k=0}^{\infty} \frac{(a, b, c ; q)_k}{(q, a q / b, a q / c ; q)_k} A_k=\sum_{j,k\geq 0}\left(\frac{a}{\lambda}\right)^j   \frac{\left(a q^{2 j}, a / \lambda, \lambda b c q^j / a ; q\right)_k}{\left(q, \lambda q^{2 j+1}, a^2 q^{j+1} / \lambda b c ; q\right)_k}\\
&\qquad\qquad\qquad\qquad\times  \frac{(\lambda ; q)_j\left(1-\lambda q^{2 j}\right)(\lambda b / a, \lambda c / a, a q / b c ; q)_j(a ; q)_{2 j}}{(q ; q)_j(1-\lambda)\left(a q / b, a q / c, q a^2 / \lambda b c ; q\right)_j(\lambda q ; q)_{2 j}} A_{j+k}.
\end{align*}
\end{description}
No doubt, these two general transformations are essentially different from L. Wang et al's result \eqref{WWWWWW}. Therefore,  it is worthwhile to investigate other transformation like \eqref{WWWWWW} in a general setting of $q$-series. As a consequence of such motivation, in this paper, by utilizing  the classical Pfaff-Saalsch$\ddot{u}$tz   summation formula of ${}_3\phi_2$ series (cf. \cite[(II.12)]{10})
\begin{align}
{ }_3 \phi_2\left[\begin{array}{c}
a, b, q^{-n} \\
c, a b c^{-1} q^{1-n}
\end{array} ; q, q\right]=\frac{(c / a, c / b ; q)_n}{(c, c / a b ; q)_n} ,\label{phi32}
\end{align}
we  will extend   \eqref{WWWWWW} to the following
\begin{dl}\label{mainthm}
Let $\left\{A_n\right\}_{n\in \mathbb{Z}}$ be any general sequence such that the infinite series below converge absolutely. Then it holds
 \begin{align}
\sum_{i, j \geq 0} \frac{A_{i-j} q^{ij} x^j}{(q ; q)_{i}(q ; q)_j}\frac{(bq^{-i}, t ; q)_{j}}{(bt; q)_{j}}=\sum_{i, j, k \geq 0} \frac{A_{i-j} q^{k(k-1)/2}(-1)^{k}x^{j+k}}{(q ; q)_i(q ; q)_j(q ; q)_k}\frac{(tq^{i+k},b;q)_j}{(bt;q)_j}.\label{WWWWWW-00-wg}
\end{align}
\end{dl}
Our paper is organized as follows. In the next section, we will show Theorem \ref{mainthm} merely
by \eqref{phi32}. Section 3 is devoted to transformation formulas for ${}_{r+1}\phi_r$ series derived from \eqref{WWWWWW-00-wg} with emphasis on those for both ${}_2\phi_1$ 
and ${}_3\phi_2$ 
series. As applications, some multi-sum   Rogers-Ramanujan-Slater identities are also presented.

Before proceeding, a bit of explanation  on notation is needed here. Throughout this paper we will use the standard notation and terminology
for $q$-series as in G. Gasper and M. Rahman
\cite{10}. As customary, the $q$-shifted factorials are
\begin{subequations}\label{notation-two-total}
\begin{align}
(a ; q)_\infty &:=\prod_{k=0}^\infty\left(1-a q^{k}\right), \label{notation-one}\\
(a ; q)_n &:=(a ; q)_{\infty} /\left(a q^n ; q\right)_{\infty} \quad (n\in \mathbb{Z}).\label{notation-two-1}
\end{align}
\end{subequations}
 In general, the multi-sum -factorial is given by
\begin{align*}
\left(a_1, a_2, \ldots, a_m ; q\right)_n:=\prod_{k=1}^m\left(a_k; q\right)_n.
\end{align*}
As usual,  the binomial coefficient $\binom{n}{k}:=n!/(k!(n-k)!)$ and its $q$-analogue  is
$$\left[n \atop k\right]_q :=\frac{(q ; q)_n}{(q ; q)_k(q ; q)_{n-k}} .
$$
Usually it is referred to as the $q$-binomial coefficient or the Gaussian binomial coefficient. Following \cite{10}, we introduce a basic hypergeometric series with the base $q: |q|<1$ and the argument $x$  as
\begin{align*}
{ }_{r} \phi_s\left[\begin{array}{c}
a_1, a_2,\ldots, a_{r} \\
b_1,b_2, \ldots, b_s
\end{array}; q, x\right]:=\sum_{n=0}^{\infty} \frac{\left(a_1,a_2,\ldots, a_{r} ; q\right)_n}{\left(b_1, b_2,\ldots, b_s ; q\right)_n}\tau(n)^{s+1-r} \frac{x^n}{(q; q)_n}.
\end{align*}
Hereafter, we write $\tau(n)$ briefly  for $(-1)^nq^{n(n-1)/2}$.

In addition, we will use the following basic relations
 (cf. \cite[(I.2) and (I.12)]{10})
 \begin{subequations}\label{notation-two-total-new}
\begin{align}
&(a;q)_{-n}=\tau(n)\frac{(q/a)^n}{(q/a;q)_n},\label{basic22}\\
&\left[n\atop k\right]_q\tau(k)=q^{nk}\frac{(q^{-n};q)_k}{(q;q)_k}\label{basic33}
\end{align}
as well as the basic relations
\begin{align}
\tau(m+n)&=\tau(m)\tau(n)q^{mn},\label{aaa}\\
(a;q)_{m+n}&=(a;q)_{m}(aq^m;q)_{n}=(a;q)_{n}(aq^n;q)_{m},\label{bbb}
\end{align}
where $ m,n\in \mathbb{Z}$.
\end{subequations}
\section{Proof of the main 
 result}
In this section,  we proceed to show our main result, i.e. Theorem \ref{mainthm}. 
\begin{pot444} It is sufficient to calculate
\begin{align*}
&\mbox{RHS of \eqref{WWWWWW-00-wg}}=\sum_{i, j\geq k\geq 0}\frac{A_{i-j}  x^j(t q^i, b;q)_{j-k}}{(q;q)_{i-k}(q, b t ; q)_{j-k}} \frac{\tau(k)}{(q;q)_k}\\
&=\sum_{i, j\geq0} \frac{A_{i-j} x^j(t q^i, b ;q)_j}
{(q;q)_i(q, b t ;q)_j}\sum_{k=0}^{\min \{i, j\}} \frac{\left(t q^{i+j}, b q^j ; q\right)_{-k} \tau(k)}{\left(q^{1+i}, q^{1+j}, b t q^j;q\right)_{-k}(q;q)_k}\\
&=\sum_{i, j \geq 0} \frac{A_{i-j} x^j(t q^i, b; q)_j}{(q;q)_i(q, b t; q)_j}{}_3 \phi_2\left[\begin{array}{l}
q^{-i}, q^{-j}, q^{1-j} / b t \\
q^{1-i-j} / t, q^{1-j} / b
\end{array};q,q\right].
\end{align*}
In this stage, we apply \eqref{bbb} to  the inner sum of the second equality and then \eqref{basic22} for the last equality. Next, by the classical Pfaff-Saalsch$\ddot{u}$tz   summation formula \eqref{phi32}, it is clear that the ${}_3\phi_2$ series on the right side
$$
{}_3 \phi_2\left[\begin{array}{l}
q^{-i}, q^{-j}, q^{1-j} / b t \\
q^{1-i-j} / t, q^{1-j} / b
\end{array};q,q\right]=\frac{(q^{1-j} / t, b q^{-i};q)_j}{(q^{1-i-j} / t, b;q)_j}.$$
As a result, it follows
\begin{align*}
\mbox{RHS of \eqref{WWWWWW-00-wg}}&=\sum_{i, j \geq 0} \frac{A_{i-j}x^j(tq^i, b; q)_j}{(q;q)_i(q, bt;q)_j}~\frac{(q^{1-j} / t, b q^{-i};q)_j}{(q^{1-i-j} / t, b;q)_j}\\
&=\sum_{i, j\geq 0} \frac{A_{i-j}x^j(b q^{-i}, t;q)_j q^{i j}}{(q;q)_ i(q, bt;q)_j}=
\mbox{LHS of \eqref{WWWWWW-00-wg}.}
\end{align*}
The theorem is proved.
\qed
\end{pot444}

Considering applications to $q$-series, it is necessary to reformulate  Theorem \ref{mainthm} in an equivalent form.
\begin{dl}
\label{mainthm-add}
Let $\{A_n\}_{n\in\mathbb{Z}}$ be any  general sequence such that the series below converge absolutely. Then  it holds
 \begin{align}
\sum_{i, j \geq 0} \frac{A_{i-j} q^{ij} x^j}{(q ; q)_{i}}\frac{(bq^{-i}, t ; q)_{j}}{(q,bt; q)_{j}}
=(x;q)_\infty\sum_{i, j,k\geq 0} \frac{A_{i-j}(tq^{i+j})^kx^{j}}{(q ; q)_i(q; q)_j(q; q)_k}
\frac{(q^{-j};q)_k(b;q)_j}{(x;q)_k(bt;q)_j}.\label{WWWWWW-00-wg-1204}
\end{align}
\end{dl}
\pf From the $q$-binomial theorem (cf. \cite[(II.3)]{10}), it follows
\[(tq^{i+k};q)_j=\sum_{l=0}^j\left[j\atop l\right]_q\tau(l) \big(tq^{i+k}\big)^l.\]
Thus we easily find
\begin{align*}
\mbox{RHS of \eqref{WWWWWW-00-wg} }&=\sum_{i, j, k \geq 0} \frac{A_{i-j} \tau(k)x^{j+k}}{(q ; q)_i(q ; q)_j(q ; q)_k}\frac{(b;q)_j}{(bt;q)_j}\bigg\{\sum_{l=0}^j\left[j\atop l\right]_q\tau(l) \big(tq^{i+k}\big)^l\bigg\}\\
&=\sum_{i, j\geq  l\geq 0} \frac{A_{i-j}\tau(l) t^lq^{il}x^{j}}{(q ; q)_i(q ; q)_l(q;q)_{j-l}}
\frac{(b;q)_j}{(bt;q)_j}\bigg\{\sum_{k=0}^\infty\frac{\tau(k)}{(q ; q)_k}(xq^l)^k\bigg\}\\
&=(x;q)_\infty\sum_{i, j\geq  l\geq 0} \frac{A_{i-j}\tau(l) (tq^{i})^lx^{j}}{(q ; q)_i(q; q)_l(q;q)_{j-l}}
\frac{(b;q)_j}{(bt;q)_j(x; q)_l}.
\end{align*}
Using \eqref{basic33} again to simplify the last sum, we have
\begin{align*}
\mbox{RHS of \eqref{WWWWWW-00-wg} }&=(x;q)_\infty\sum_{i, j,l\geq 0} \frac{A_{i-j}(tq^{i+j})^lx^{j}}{(q ; q)_i(q; q)_j(q; q)_l}
\frac{(q^{-j};q)_l(b;q)_j}{(bt;q)_j(x;q)_l}.
\end{align*}
Replacing $l$ with $k$ leads us to \eqref{WWWWWW-00-wg-1204}.
\qed

\section{Applications}
This section is devoted to  applications of Theorems \ref{mainthm} and \ref{mainthm-add} to $q$-series.
\subsection{Various specific transformations and examples}
Now we proceed to investigate some special cases of Theorems \ref{mainthm} and  \ref{mainthm-add}, as well as their potential  applications.
\subsubsection{Various specific transformations}
 First of all, we need to show Theorem  \ref{mainthm}
  indeed contains Theorem \ref{wanngliuquan} as a special case.
 \begin{pot333}
 It suffices to take the limit as $t\to \infty$ first and then $b\to \infty$ in \eqref{WWWWWW-00-wg}. Note that
\[\lim_{t\to\infty}\frac{(tq^{i+k},b;q)_j}{(bt;q)_j}=\frac{(b;q)_j}{b^j}q^{(i+k)j}
\quad\mbox{and}~~\lim_{b\to\infty}\frac{(b;q)_j}{b^j}=\tau(j).\]
All these reduces the right-hand side of \eqref{WWWWWW-00-wg} to
\begin{align*}\sum_{i, j, k \geq 0} \frac{A_{i-j} (-x)^{j+k}}{(q ; q)_i(q ; q)_j(q ; q)_k}q^{k(k-1)/2+(i+k)j+j(j-1)/2}\\=(x;q)_{\infty}\sum_{i, j\geq 0} \frac{A_{i-j}q^{ij+j(j-1)/2} (-x)^{j}}{(q ; q)_i(q ; q)_j(x ; q)_j}\end{align*}
and the left-hand side of \eqref{WWWWWW-00-wg} to
\[\sum_{i, j \geq 0} \frac{A_{i-j} q^{j(j-1)/2} (-x)^j}{(q ; q)_{i}(q ; q)_j}.\]
As a final step, replace $x$ with $-x$. The proof is finished.
\qed
\end{pot333}

The special case $t=1$ of  \eqref{WWWWWW-00-wg} in
  Theorem \ref{mainthm} deserves our particular attention.
\begin{tl}Under the same conditions as Theorem \ref{mainthm}. We have
\begin{align}\sum_{n=0}^\infty \frac{A_{n} }{(q;q)_{n}}=\sum_{i, j,k\geq 0} \frac{A_{i-j}(q^{i+k};q)_j}{(q ; q)_i(q; q)_j(q; q)_k}\tau(k)x^{j+k}.\label{WWWWWW-00000-cw}
\end{align}
\end{tl}

From Theorem \ref{mainthm-add} it is easy to deduce a few easy-to-use transformations.  For instance,  by taking the limit as  $t\to\infty$ alone in \eqref{WWWWWW-00-wg-1204}, we have
\begin{tl} \label{mainthm-add-add}With the same condition as Theorem \ref{mainthm-add}. Then
 it holds
 \begin{align}
\sum_{i, j \geq 0} \frac{A_{i-j} q^{ij} (x/b)^j}{(q ; q)_{i}}\frac{(bq^{-i}; q)_{j}}{(q; q)_{j}}
=(x;q)_\infty\sum_{i, j\geq 0} \frac{A_{i-j}q^{ij}(x/b)^{j}}{(q ; q)_i(q; q)_j}
\frac{(b;q)_j}{(x;q)_j}.\label{WWWWWW-00-wg-1204-1}
\end{align}
\end{tl}
Transformation \eqref{WWWWWW-00-wg-1204-1} may be further stated as
 \begin{tl}Let $\left\{A_n\right\}_{n\in \mathbb{Z}}$ be any general sequence such that the series below converge absolutely. Then there holds
\begin{align}
&\sum_{i \geq 0} \frac{(q/b;q)_i y^i}{(q ; q)_{i}}\bigg(\sum_{j \geq 0}\frac{A_{i-j}\tau(j)(x/y)^{j}}{(q; q)_{j}}\bigg)\nonumber\\
&=(x;q)_\infty\sum_{j\geq 0}
\frac{\tau(j)(x/y)^{j}}{(q; q)_j(x;q)_j}\bigg(\sum_{i\geq 0}\frac{A_{i-j}(q^{1-j}/b;q)_i(q^jy)^{i}}{(q ; q)_i}\bigg).\label{WWWWWW-0000000-1-2}
\end{align}
\end{tl}
\pf It suffices to replace
$A_i$  with $(q/b;q)_iy^iA_i$ in \eqref{WWWWWW-00-wg-1204-1}.  In this case,  it is easy to check
\begin{align}(q/b;q)_{i-j}(bq^{-i};q)_j&=b^j\tau(j)q^{-ij}(q/b;q)_i,\label{xxxrma-one}\\
(q/b;q)_{i-j}(b;q)_j&=b^j\tau(j)(q^{1-j}/b;q)_i.\nonumber
\end{align}
Thus we may reformulate \eqref{WWWWWW-00-wg-1204-1} as
\begin{align*}
\sum_{i, j \geq 0} \frac{(q/b;q)_i y^i}{(q ; q)_{i}}\frac{A_{i-j}\tau(j)(x/y)^{j}}{(q; q)_{j}}
&=(x;q)_\infty\sum_{i, j\geq 0} \frac{\tau(j)(q^{1-j}/b;q)_iy^{i}}{(q ; q)_i}
\frac{(x/y)^{j}A_{i-j}q^{ij}}{(q; q)_j(x;q)_j},
\end{align*}
which is in agreement with  \eqref{WWWWWW-0000000-1-2}.
\qed

It is also of interest that letting $b\to \infty$ in \eqref{WWWWWW-0000000-1-2}, we recover \eqref{WWWWWW} up to the factor $y^{i-j}$. Once certain additional conditions are satisfied by the sequence $\{A_{n}\}_{n\in\mathbb{Z}}$, we may find some useful transformations.  To justify this, we  now  present   a convolutional-form transformation using  \eqref{WWWWWW-00-wg-1204-1}.
\begin{tl}Assume that $A_{n}=0$ for $n<0$. Then we have
\begin{align}\sum_{n=0}^\infty\frac{(q/b;q)_n}{(q ; q)_{n}}\boldsymbol\lbrack t^n\boldsymbol\rbrack\bigg\{(xt;q)_\infty\sum_{i=0}^\infty \frac{A_{i}\tau(i)}{(q/b; q)_{i}}(tq)^i \bigg\}\nonumber\\
=(x;q)_\infty\sum_{n=0}^\infty\frac{\tau(-n)}{(q ; q)_n}\boldsymbol\lbrack t^n\boldsymbol\rbrack\bigg\{{}_{1}\phi_1\left[{b\atop x}; q, \frac{xt}{b}\right]\sum_{i=0}^\infty A_{i}t^i\bigg\},\label{WWWWWW-00-wg-1204-1-newidea}
\end{align}
where the coefficient functional $\boldsymbol\lbrack t^n\boldsymbol\rbrack$ over the ring of Laurent power series $\mathbb{C}[[t]]$ is defined by
\[\boldsymbol\lbrack t^n\boldsymbol\rbrack\sum_{k=-\infty}^\infty a_kt^k=a_n.\]
\end{tl}
\pf Clearly, using the relation \eqref{xxxrma-one},
we may rewrite
\begin{align*}
\mbox{LHS of   \eqref{WWWWWW-00-wg-1204-1}}&=
\sum_{i\geq 0} \frac{(q/b;q)_i}{(q ; q)_{i}}\sum_{j =0}^i \frac{A_{i-j}}{(q/b; q)_{i-j}} \times \frac{\tau(j)x^j}{(q; q)_{j}}\\
&= \sum_{n\geq 0} \frac{(q/b;q)_n}{(q ; q)_{n}}\boldsymbol\lbrack t^n\boldsymbol\rbrack\bigg\{(xt;q)_\infty\sum_{i=0}^\infty \frac{A_{i}}{(q/b; q)_{i}}t^i \bigg\}.
\end{align*}
On the other hand, we easily check that
\begin{align*}\mbox{RHS of   \eqref{WWWWWW-00-wg-1204-1}}&=(x;q)_\infty\sum_{i\geq 0} \frac{1}{(q ; q)_i} \sum_{j=0}^i\frac{A_{i-j}\tau(j-i)q^{ij}}{\tau(j-i)}\frac{(x/b)^{j}(b;q)_j}{(q,x;q)_j} \\
&=(x;q)_\infty\sum_{i\geq 0} \frac{\tau(-i)}{(q ; q)_i} \sum_{j=0}^i\frac{A_{i-j}}{\tau(j-i)}\frac{\tau(j)(x/b)^{j}(b;q)_j}{(q,x;q)_j} \\
&=(x;q)_\infty\sum_{n\geq 0} \frac{\tau(-n)}{(q ; q)_n}\boldsymbol\lbrack t^n\boldsymbol\rbrack\bigg\{{}_{1}\phi_1\bigg[{b \atop x}; q, \frac{xt}{b}\bigg]\sum_{i=0}^\infty\frac{A_{i}t^i}{\tau(-i)}\bigg\}.
\end{align*}
In order to ensure the convergency of the sums related, we replace $A_i$ with $A(i)\tau(-i)$. Thus we get the desired transformation \eqref{WWWWWW-00-wg-1204-1-newidea}.

Alternatingly, let $t=1$ in \eqref{WWWWWW-00-wg-1204}. We can find another two new transformations similar to \eqref{WWWWWW-00000-cw}.
\begin{tl}\label{secondtrans}
Let $\left\{A_n\right\}_{n\in \mathbb{Z}}$  be any sequence such that the series below converge absolutely. Then there hold
\begin{align}\sum_{n=0}^\infty \frac{A_{n} }{(q ; q)_{n}}
=(x;q)_\infty\sum_{i,k\geq 0} \frac{x^kq^{ik}\tau(k)}{(q ; q)_i(q,x; q)_k}\sum_{ j\geq 0}
\frac{A_{i-k-j}}{(q;q)_{j}}x^{j},\label{ssssss-1}\\
\sum_{n=0}^\infty \frac{A_{-n} }{(q ; q)_{n}}
=(x;q)_\infty\sum_{i,k\geq 0} \frac{x^kq^{ik}\tau(k)}{(q ; q)_i(q,x; q)_k}\sum_{j\geq 0}
\frac{A_{j-(i-k)}}{(q;q)_{j}}x^{j}.\label{ssssss-1-added}
\end{align}
\end{tl}
\pf Actually, when $t=1$,   \eqref{WWWWWW-00-wg-1204} reduces to
 \begin{align}
\sum_{n=0}^\infty \frac{A_{n} }{(q ; q)_{n}}
&=(x;q)_\infty\sum_{i, j,k\geq 0} \frac{A_{i-j}q^{ik}x^{j}}{(q ; q)_i(q; q)_j(x; q)_k}
\frac{(q^{-j};q)_kq^{jk}}{(q;q)_k}\label{128}\\
&=(x;q)_\infty\sum_{i, j,k\geq 0} \frac{A_{i-j}q^{ik}x^{j}}{(q ; q)_i(x; q)_k}
\frac{\tau(k)}{(q;q)_{j-k}(q;q)_k}\nonumber\\
&=(x;q)_\infty\sum_{i,k\geq 0} \frac{x^kq^{ik}\tau(k)}{(q ; q)_i(x; q)_k(q;q)_k}\sum_{J:= j-k\geq 0}
\frac{A_{i-k-J}x^{J}}{(q;q)_{J}}.\nonumber
\end{align}
Thus \eqref{ssssss-1} is proved.
Evidently, since $\{A_i\}_{i\in\mathbb{Z}}$ is arbitrary, we can replace $A_i$ with  $A_{-i}$, thus obtaining  \eqref{ssssss-1-added}.
\qed

Furthermore, the special case $b=t=0$ of
\eqref{WWWWWW-00-wg-1204}  of Theorem \ref{mainthm-add} yields
\begin{tl}Under the same conditions as Theorem \ref{mainthm-add}. We have
\begin{align}
\sum_{i, j \geq 0} \frac{A_{i-j} q^{ij} x^j}{(q ; q)_{i}(q ; q)_j}=(x;q)_\infty\sum_{i, j\geq 0} \frac{A_{i-j}x^{j}}{(q ; q)_i(q ; q)_j}.\label{WWWWWW-000}
\end{align}
\end{tl}

\subsubsection{Double and triple Rogers-Ramanujan-Slater identities}
In the part, we will consider  some more specific  results via specializing $\{A_n\}_{n\in \mathbb{Z}}$. First, there are two concrete identities which are worthy of separate corollaries.
\begin{lz} Set $A_i=\tau(i)(x/y)^i$ in \eqref{ssssss-1}. Then the Cauchy identity follows
\begin{align}\sum_{n=0}^\infty \frac{q^{n^2}x^n}{(q,xq;q)_{n}}=\frac{1}{(xq;q)_\infty}
.\label{ssssss-199}
\end{align}
Instead, take $A_i=1/y^i$ in \eqref{ssssss-1}. Then we get
\begin{align*}\sum_{n=0}^\infty \frac{(1/y;q)_n}{(q,x; q)_n}\tau(n)(xy)^n=\frac{(xy;q)_\infty}{(x;q)_\infty}
.
\end{align*}
It is not only a generalization of the Cauchy identity \eqref{ssssss-199} but also the limit $b\to\infty$ of the Gauss summation formula  \cite[(II.8)]{10}, i.e.
\begin{align*}\sum_{n=0}^\infty \frac{(1/y,b;q)_n}{(q,x; q)_n}\bigg(\frac{xy}{b}\bigg)^n=\frac{(xy,x/b;q)_\infty}{(x,xy/b;q)_\infty}
.
\end{align*}
\end{lz}
Based on \eqref{WWWWWW-000},  it is easy to deduce a few results on multi-sum   Rogers-Ramanujan-Slater type identities in \cite{12-jmaa} given by Z. Cao and L. Wang.
\begin{lz}[cf. \mbox{\cite[Thm. 3.1]{12-jmaa}}]
\begin{align}\sum_{i,j\geq 0} \frac{\tau(i-j)y^{i-j}x^j}{(q ; q)_i(q; q)_j}=\frac{(y;q)_\infty(xq/y;q)_\infty}{(x;q)_\infty}.\label{WWWWWW-0000000-1}
\end{align}
\end{lz}
\pf It only needs  to set
$A_i=\tau(i)y^{i}$ in \eqref{WWWWWW-000}.
\qed

 \begin{lz}[cf. \mbox{\cite[Thm. 3.2]{12-jmaa}}]
\begin{align}\sum_{i, j\geq 0} \frac{(-1)^{i-j}q^{(i-j)^2}x^{j}}{(q^2 ; q^2)_i(q^2 ; q^2)_j}=\frac{(q;q^2)_\infty(xq;q^2)_\infty}{(x;q^2)_\infty}.\label{WWWWWW-0000000-111}
\end{align}
\end{lz}
\pf To achieve \eqref{WWWWWW-0000000-111}, we first need  to  replace $q$ with $q^2$, and then take $A_i=(-1)^{i}q^{i^2}$ in  \eqref{WWWWWW-000}.
So we have
\begin{align*}\sum_{i, j\geq 0} \frac{(-1)^{i-j}q^{(i-j)^2}x^{j}}{(q^2 ; q^2)_i(q^2 ; q^2)_j}&=\frac{1}{(x;q^2)_\infty}
\sum_{i,j \geq 0} \frac{(-1)^{i-j}q^{i^2+j^2}  x^j}{(q^2 ; q^2)_{i}(q^2 ; q^2)_j}\\
&=\frac{1}{(x;q^2)_\infty}
\sum_{i \geq 0} \frac{(-1)^{i}q^{i^2}}{(q^2 ; q^2)_{i}}
\sum_{j \geq 0} \frac{(-1)^{j}q^{j^2}  x^j}{(q^2 ; q^2)_j}=\frac{(q;q^2)_\infty(xq;q^2)_\infty}{(x;q^2)_\infty}.
\end{align*}
As desired.
\qed

Moreover, we can deduce Theorems 3.4 and 3.6 of \cite{12-jmaa} from our \eqref{WWWWWW-000}. To do this, we need to show
\begin{yl}For any complex numbers $x,y: |x|<1$, it holds
\begin{align}
\sum_{n=0}^\infty\frac{x^n}{(q,y;q)_n}=\frac{1}{(x,y;q)_\infty}\sum_{n=0}^\infty\frac{(x;q)_n}{(q;q)_n}\tau(n)y^n.
\label{addedidentity}
\end{align}
In particular, the following limits hold:
\begin{subequations}\label{ffffff}
\begin{align}
\lim_{x\to 1}(x;q)_\infty\sum_{n=0}^\infty\frac{x^n}{(q,y;q)_n}&=\frac{1}{(y;q)_\infty},\label{addedidentity-1}\\
\lim_{x\to q}(x;q)_\infty\sum_{n=0}^\infty\frac{x^n}{(q,y;q)_n}&=\frac{1}{(y;q)_\infty}\sum_{n\geq 0}\tau(n)y^n\label{addedidentity-q},\\
\lim_{x\to 1/q}(x;q)_\infty\sum_{n=0}^\infty\frac{x^n}{(q,y;q)_n}&=\frac{1+y/q}{(y;q)_\infty}.\label{addedidentity-nq}
\end{align}
\end{subequations}
\end{yl}
\pf To establish \eqref{addedidentity}, it suffices to show
\begin{align*}
\frac{1}{(q,y;q)_m}=\frac{1}{(y;q)_\infty}\sum_{n=0}^\infty\frac{\tau(n)y^n}{(q;q)_n}\frac{q^{mn}}{(q;q)_m},
\end{align*}
which follows by  equating the coefficients of the term $x^m$  on both sides of \eqref{addedidentity}.
A simplification leads us to the nonterminating  $q$-binomial theorem:
\begin{align*}
(yq^m;q)_\infty=\sum_{n=0}^\infty\frac{\tau(n)}{(q;q)_n}(yq^m)^{n}.
\end{align*}
It is obvious that when $x=1$ or $x=1/q$, the corresponding sum on the left-hand side of \eqref{addedidentity} is nonconvergent, but the limit as $x\to a~~(a=1,1/q)$ of \eqref{ffffff} do exist, as we wanted.
\qed

Now we are able to recover successfully
\begin{lz}[cf. \mbox{\cite[Thm. 3.4]{12-jmaa}}]\label{zzzzzz}Define $\tau_p(n):=(-1)^nq^{p\binom{n}{2}}$. Then it holds
\begin{align}\sum_{i,j\geq 0} \frac{\tau_p(j-i)\tau(i-j)y^{i-j}q^{ij}}{(q ; q)_i(q; q)_j}
=\frac{(-q/y,-q^py,q^{p+1};q^{p+1})_\infty}{(q;q)_\infty}.\label{WWWWWW-0000000-11}
\end{align}
\end{lz}
\pf It suffices to take
$A_i=\tau_p(-i)\tau(i)y^{i}$ in  \eqref{WWWWWW-000}. Then it becomes
\begin{align}
\sum_{i, j\geq 0} \frac{\tau_p(j-i)\tau(i-j)y^{i-j} q^{ij}x^{j}}{(q ; q)_i(q ; q)_j}=(x;q)_\infty
\sum_{i, j \geq 0} \frac{\tau_p(j-i)\tau(i-j)y^{i-j} x^{j} }{(q ; q)_{i}(q ; q)_j}\label{GGGGGG}.
\end{align}
Hence, a bit simplification yields
\begin{align*}
\mbox{RHS of \eqref{GGGGGG}}&=(x;q)_\infty
\sum_{j-i=n; i,j \geq  0} \frac{\tau_p(n)\tau(-n)y^{-n} x^{j} }{(q ; q)_{i}(q ; q)_j}=S_1(x)+S_2(x),
\end{align*}
where
\begin{align*}
S_1(x)&:=(x;q)_\infty
\sum_{j-i=n\geq 0; i,j \geq  0} \frac{\tau_p(n)\tau(-n)y^{-n} x^{j} }{(q ; q)_{i}(q ; q)_{j}},\\
S_2(x)&:=(x;q)_\infty
\sum_{j-i=n\leq -1; i,j \geq  0} \frac{\tau_p(n)\tau(-n)y^{-n} x^{j} }{(q ; q)_{i}(q ; q)_j}.
\end{align*}
Next, it is easy to compute
\begin{align*}
S_1(x)&=\sum_{n\geq 0} \frac{\tau_{p+1}(n)(-xq/y)^{n} }{(q;q)_n}(x;q)_\infty\sum_{i\geq 0}\frac{x^{i} }{(q ; q)_{i}(q^{n+1} ; q)_{i}},\\
S_2(x)&=\sum_{n\geq 1} \frac{\tau_{p+1}(-n)(-y/q)^{n} }{(q;q)_{n}}(x;q)_\infty\sum_{j\geq 0}\frac{x^{j} }{(q ; q)_{j}(q^{n+1} ; q)_{j}}.
\end{align*}
As a final step, take all these expressions into account and let $x\to 1$ on both sides of \eqref{GGGGGG}. According to \eqref{addedidentity-1}, we have
\begin{align*}\lim_{x\to 1}S_1(x)&=\sum_{n\geq 0} \frac{\tau_{p+1}(n)(-q/y)^{n} }{(q;q)_n}\frac{1}{(q^{n+1};q)_\infty}=\frac{1}{(q;q)_\infty}\sum_{n\geq 0}\tau_{p+1}(n)(-q/y)^{n},\\
\lim_{x\to 1}S_2(x)&=\sum_{n\geq 1} \frac{\tau_{p+1}(-n)(-y/q)^{n} }{(q;q)_n}\frac{1}{(q^{n+1};q)_\infty}=\frac{1}{(q;q)_\infty}\sum_{n\geq 1}\tau_{p+1}(-n)(-y/q)^{n},
\end{align*}
obtaining
\begin{align*}
\mbox{LHS of \eqref{WWWWWW-0000000-11}}=\frac{1}{(q;q)_\infty}\sum_{n\in\mathbb{Z}}\tau_{p+1}(n)(-q/y)^{n}=
\mbox{RHS of \eqref{WWWWWW-0000000-11}}.
\end{align*}
Note that the  second equality comes from Jacobi's triple product identity \cite[(II. 28)]{10}, namely
\begin{align}
\sum_{n\in\mathbb{Z}}\tau(n)x^{n}=(x,q/x,q;q)_\infty.\label{jacpbintriple}
\end{align}
\qed

As mentioned earlier, we can recover Theorem 3.6 of \cite{12-jmaa}  by the same method.
\begin{lz}[cf. \mbox{\cite[Thm. 3.6]{12-jmaa}}]\label{xxxxxx}With the same notation as above. Then there hold
\begin{align}&\sum_{i,j\geq 0} \frac{\tau_p(j-i)\tau(i-j)y^{i-j}q^{(i-1)j}}{(q ; q)_i(q; q)_j}\label{WWWWWW-0000000-1124}\\
&\qquad=
\frac{(-1/y,-q^{p+1}y,q^{p+1};q^{p+1})_\infty}{(q;q)_\infty}+
\frac{(-q/y,-q^py,q^{p+1};q^{p+1})_\infty}{(q;q)_\infty}.\nonumber
\end{align}
\end{lz}
\pf Actually, \eqref{WWWWWW-0000000-1124} comes from the limits as $x\to 1/q$ on both sides of \eqref{GGGGGG}, stated by \eqref{addedidentity-nq}. Since, in such case, it follows
\begin{align*}\mbox{LHS of \eqref{WWWWWW-0000000-1124}
}
&=\frac{1}{(q;q)_\infty}\bigg(\sum_{n\geq 0}\tau_{p+1}(n)(-1/y)^{n}(1+q^n)+\sum_{n\geq 1}\tau_{p+1}(-n)(-y/q)^{n} (1+q^n)\bigg)\\
&=\frac{1}{(q;q)_\infty}\bigg(\sum_{n\geq 0}\tau_{p+1}(n)(-1/y)^{n}+\sum_{n\geq 1}\tau_{p+1}(-n)(-y/q)^{n} \nonumber\\
&\qquad\qquad~~+\sum_{n\geq 0}\tau_{p+1}(n)(-q/y)^{n}+\sum_{n\geq 1}\tau_{p+1}(-n)(-y)^{n}\bigg)\nonumber\\
&=\frac{1}{(q;q)_\infty}\bigg(\sum_{n\in\mathbb{Z}}\tau_{p+1}(n)(-1/y)^{n}+\sum_{n\in\mathbb{Z}}\tau_{p+1}(n)(-q/y)^{n}\bigg)\nonumber.
\end{align*}
By applying Jacobi's triple product identity \eqref{jacpbintriple} to both bilateral series, we obtain \eqref{WWWWWW-0000000-1124} at once.
\qed

In the same vein  as Examples \ref{zzzzzz} and  \ref{xxxxxx}, we have the following new multi-sum   Rogers-Ramanujan-Slater identity.
\begin{lz}With the same notation as above. Then it holds
\begin{align}&\sum_{i,j\geq 0} \frac{\tau_p(j-i)\tau(i-j)y^{i-j}q^{(i+1)j}}{(q ; q)_i(q; q)_j}\label{WWWWWW-0000000-1123}\\
&\qquad=\frac{1}{(q;q)_\infty}\bigg(\bigg\{\sum_{n\geq 0}\tau_{p+1}(n)(-q^2/y)^{n}+\sum_{n\geq 1}\tau_{p+1}(-n)(-y/q)^{n} \bigg\}\sum_{i\geq 0}\tau(i) q^{(n+1)i}\bigg)\nonumber
\end{align}
\end{lz}
\pf It comes from the limits as $x\to q$ on both sides of \eqref{GGGGGG}, stated by \eqref{addedidentity-q}. 
\qed

Actually, we are always able to choose such $A_n$ that the left-hand side of \eqref{128} takes a closed form. Here are two basic examples.
\begin{lz}For complex numbers $a,b,c: |c/ab|<1$, it holds
\begin{align}\sum_{i, j, k\geq 0}\frac{(a,b;q)_{i}}{(q,c; q)_i} \frac{(q^{1-i}/c;q)_{j}(q^{-j};q)_k}{(q,q^{1-i}/a,q^{1-i}/b;q)_{j}}\frac{q^{ik+jk-\binom{j}{2}}}
{(q,x; q)_k}\bigg(\frac{c}{ab}\bigg)^i (-x)^j=
\frac{(c/a,c/b;q)_\infty}{(x,c,c/ab;q)_\infty}.\label{WWWWWW-0000000-new}
\end{align}
\end{lz}
\pf It suffices to take in \eqref{128} that
$$A_i=\frac{(a,b;q)_i}{(c;q)_i}\bigg(\frac{c}{ab}\bigg)^i.
$$
And then \eqref{WWWWWW-0000000-new} comes from  the Gauss summation formula (cf. \cite[(II.8)]{10})
\begin{align}
\sum_{n=0}^\infty \frac{(a,b;q)_n}{(q,c;q)_n}\bigg(\frac{c}{ab}\bigg)^n=
\frac{(c/a,c/b;q)_\infty}{(c,c/ab;q)_\infty}.\label{gauss-id}
\end{align}
\qed

\begin{lz}
\begin{align}\sum_{i, j, k\geq 0} \frac{(a;q)_i(q^{-j};q)_kq^{ik+jk+i^2/2+i/2+j^2-2ij}}{(q ; q)_i(q,q^{1-i}/a; q)_j(q,x; q)_k}\bigg(-\frac{x}{a}\bigg)^j=\frac{(-q;q)_\infty(aq;q^2)_\infty}{(x;q)_\infty}.\label{WWWWWW-0000000}
\end{align}
\end{lz}
\pf It suffices to take in \eqref{128}
$A_i=    q^{i(i+1)/2}(a;q)_i.
$
And then \eqref{WWWWWW-0000000} comes from  the Lebesgue identity (cf. \cite[Ex. 1.16]{10})
$$\sum_{n=0}^\infty \frac{q^{n(n+1)/2}(a;q)_n}{(q ; q)_{n}}=(-q;q)_\infty(aq;q^2)_\infty.$$
\qed

Even more interesting to us is  the case $A_n=\delta_n$ (the Kronecker delta) of \eqref{WWWWWW-00-wg} yields
\begin{lz}
 \begin{align}\sum_{n=0}^\infty \frac{q^{n^2} x^n}{(q ; q)_{n}^2}\frac{(bq^{-n}, t ; q)_{n}}{(bt; q)_{n}}=\sum_{j, k \geq 0} \frac{\tau(k)x^{j+k}}{(q ; q)_j^2(q ; q)_k}\frac{(tq^{j+k},b;q)_j}{(bt;q)_j}.\label{WWWWWW-07}
\end{align}
\end{lz}
The case $bt=q$  and $-q$ of \eqref{WWWWWW-07}, respectively,  gives
\begin{lz}
 \begin{align}\sum_{n=0}^\infty \frac{(t ; q)_{n}^2}{(q ; q)_{n}^3}\bigg(\frac{xq}{t}\bigg)^n\tau(n)&=\sum_{j, k \geq 0} \frac{(tq^{j+k},q/t;q)_j}{(q ; q)_j^3(q ; q)_k}\tau(k)x^{j+k},\label{WWWWWW-08}\\\sum_{n=0}^\infty \frac{\tau(n)(bx)^n}{(q ; q)_{n}}\frac{(q^2/b^2; q^2)_{n}}{(q^2; q^2)_{n}}&=\sum_{j, k \geq 0} \frac{\tau(k)x^{j+k}}{(q ; q)_j(q ; q)_k}\frac{(-q^{j+k+1}/b,b;q)_j}{(q^2;q^2)_j}.\label{WWWWWW-09}
\end{align}
\end{lz}
\subsection{Transformations between ${}_{r+1}\phi_r$ series and applications}
In this subsection, we will consider Theorem \ref{mainthm} specialized by the sequence
\begin{align}A_{n}=\frac{(a_1,a_2,\cdots,a_{r+1};q)_n}{(b_1,b_2,\cdots,b_r;q)_n}y^n.\label{rrrrrr}
\end{align}

Indeed,  replace  $A_n$ in  \eqref{WWWWWW-00-wg} with  \eqref{rrrrrr}. After a bit simplification via the  basic relations \eqref{notation-two-total-new}, it is easily seen
\begin{dl}\label{firsttrans}
\begin{align}
&\sum_{n=0}^\infty \frac{\tau(n)}{(q ; q)_n}\frac{(b,t,q/b_1,q/b_2,\cdots,q/b_r;q)_n}{(bt,q/a_1,q/a_2,\cdots,q/a_{r+1};q)_n}
\bigg(\frac{xqb_1b_2\cdots b_r}{ya_1a_2\cdots a_{r+1}}\bigg)^n\nonumber\\
&\qquad\qquad\times{ }_{r+2} \phi_{r+1}\left[\begin{array}{c}
a_1q^{-n}, \ldots, a_{r+1}q^{-n},q/b \\
b_1q^{-n}, \ldots, b_rq^{-n},q^{1-n}/b
\end{array}; q, y\right]\nonumber\\
&=\sum_{k \geq 0} \frac{\tau(k)x^{k}}{(q ; q)_k}\sum_{n=0}^\infty \frac{\tau(n)}{(q; q)_n}\frac{(q/b_1,q/b_2,\cdots,q/b_r,b,tq^k;q)_n}{(q/a_1,q/a_2,\cdots,q/a_{r+1},bt;q)_n}
\bigg(\frac{xqb_1b_2\cdots b_r}{ya_1a_2\cdots a_{r+1}}\bigg)^n\nonumber\\
&\qquad\qquad\qquad\times{ }_{r+2} \phi_{r+1}\left[\begin{array}{c}
a_1q^{-n}, \ldots, a_{r+1}q^{-n},tq^{k+n} \\
b_1q^{-n}, \ldots, b_rq^{-n},tq^{k}
\end{array}; q, y\right].\label{WWWWWW-999}
\end{align}
In particular,
let $t\to \infty$ alone, we have
\begin{align}
&\sum_{n=0}^\infty \frac{\tau(n)}{(q ; q)_n}\frac{(b,q/b_1,q/b_2,\cdots,q/b_r;q)_n}{(q/a_1,q/a_2,\cdots,q/a_{r+1};q)_n}
\bigg(\frac{xqb_1b_2\cdots b_r}{yba_1a_2\cdots a_{r+1}}\bigg)^n\nonumber\\
&\qquad\qquad\times{ }_{r+2} \phi_{r+1}\left[\begin{array}{c}
a_1q^{-n}, \ldots, a_{r+1}q^{-n},q/b \\
b_1q^{-n}, \ldots, b_rq^{-n},q^{1-n}/b
\end{array}; q, y\right]\nonumber\\
&=(x ; q)_{\infty} \sum_{n=0}^\infty \frac{\tau(n)}{(q,x ; q)_n} \frac{(q/b_1,q/b_2,\cdots,q/b_r,b;q)_n}{(q/a_1,q/a_2,\cdots,q/a_{r+1};q)_n}
\bigg(\frac{xqb_1b_2\cdots b_r}{yba_1a_2\cdots a_{r+1}}\bigg)^n\nonumber\\
&\qquad\qquad\qquad\times{ }_{r+1} \phi_{r}\left[\begin{array}{c}
a_1q^{-n}, \ldots, a_{r+1}q^{-n}\\
b_1q^{-n}, \ldots, b_rq^{-n}
\end{array}; q, yq^n\right]\label{WWWWWW-99}
\end{align}
while $t\to \infty$ and then $b\to \infty$, we have
\begin{align}
&\sum_{n=0}^\infty \frac{q^{n^2}}{(q ; q)_n}\frac{(q/b_1,q/b_2,\cdots,q/b_r;q)_n}{(q/a_1,q/a_2,\cdots,q/a_{r+1};q)_n}
\bigg(\frac{xb_1b_2\cdots b_r}{ya_1a_2\cdots a_{r+1}}\bigg)^n\nonumber\\
&\qquad\qquad\times{ }_{r+1} \phi_r\left[\begin{array}{c}
a_1q^{-n}, \ldots, a_{r+1}q^{-n} \\
b_1q^{-n}, \ldots, b_rq^{-n}
\end{array}; q, y\right]\nonumber\\
&=(x ; q)_{\infty} \sum_{n=0}^\infty \frac{q^{n^2}}{(q,x ; q)_n}\frac{(q/b_1,q/b_2,\cdots,q/b_r;q)_n}{(q/a_1,q/a_2,\cdots,q/a_{r+1};q)_n}
\bigg(\frac{xb_1b_2\cdots b_r}{ya_1a_2\cdots a_{r+1}}\bigg)^n\nonumber\\
&\qquad\qquad\qquad\times{ }_{r+1} \phi_r\left[\begin{array}{c}
a_1q^{-n}, \ldots, a_{r+1}q^{-n} \\
b_1q^{-n}, \ldots, b_rq^{-n}
\end{array}; q, yq^n\right].\label{WWWWWW-9}
\end{align}

\end{dl}
\pf All these can be verified by the basic relations \eqref{notation-two-total-new} and routine computation. The details are omitted and left to the reader.
\qed

Similarly, from \eqref{ssssss-1} it follows
\begin{dl}
 \begin{align}
&{ }_{r+1} \phi_r\left[\begin{array}{c}
a_1, \ldots, a_{r+1} \\
b_1, \ldots, b_r
\end{array}; q, y\right]
\label{ssssss}\\
&=(x;q)_\infty\sum_{i,k\geq 0} \frac{x^kq^{ik}\tau(k)}{(q ,q,x; q)_k} \frac{(a_1,a_2,\cdots,a_{r+1};q)_{i-k}}{(q^{k+1},b_1,b_2,\cdots,b_r;q)_{i-k}}y^{i-k}\nonumber\\
&\qquad\qquad\times{ }_{r+1} \phi_{r+1}\left[\begin{array}{c}0,
q^{1+k-i}/b_1,q^{1+k-i}/b_2,\cdots,q^{1+k-i}/b_r\\
q^{1+k-i}/a_1,q^{1+k-i}/a_2,\cdots,q^{1+k-i}/a_{r+1}
\end{array}; q, \frac{xb_1b_2\cdots b_rq}{ya_1a_2\cdots a_{r+1}}\right].\nonumber
\end{align}
\end{dl}
\pf  It only needs to substitute $A_n$ given by \eqref{rrrrrr} into \eqref{ssssss-1}. \qed

By replacing $A_n$ of \eqref{ssssss-1-added} with \eqref{rrrrrr},  we also easily find

\begin{dl}
 \begin{align}
 &{ }_{r+1} \phi_{r+1}\left[\begin{array}{c}0,
q/b_1,q/b_2,\cdots,q/b_r\\
q/a_1,q/a_2,\cdots,q/a_{r+1}
\end{array}; q, \frac{qb_1b_2\cdots b_r}{ya_1a_2\cdots a_{r+1}}\right]\label{ssssss-added}\\
&=(x;q)_\infty\sum_{i,k\geq 0} \frac{x^kq^{ik}\tau(k)}{(q;q)_i(q,x; q)_k}
\frac{(a_1,a_2,\cdots,a_{r+1};q)_{k-i}}{(b_1,b_2,\cdots,b_r;q)_{k-i}} y^{k-i}\nonumber\\
&\qquad\qquad\times{ }_{r+1} \phi_{r}\left[\begin{array}{c}
q^{k-i}a_1,q^{k-i}a_2,\cdots,q^{k-i}a_{r+1}\\
q^{k-i}b_1,q^{k-i}b_2,\cdots,q^{k-i}b_r
\end{array}; q,xy\right].\nonumber
\end{align}
\end{dl}
\pf  It only needs to substitute $A_n$ given by \eqref{rrrrrr} into \eqref{ssssss-1-added}.
\qed

In the rest of this paper, we will focus on possible applications of Theorem \ref{firsttrans} to $q$-series. First of all, from \eqref{WWWWWW-9} it follows a very useful transformation.
\subsubsection{Transformations related to ${}_{1}\phi_0$ series}
\begin{tl}
\begin{align}
\sum_{n=0}^\infty\frac{(z/y;q)_n}{(q,z;q)_n}\tau(n)
x^n=(x;q)_\infty \sum_{n=0}^\infty \frac{\tau(n)^2(y; q)_n}{(q,x,z; q)_n}
\bigg(\frac{xz}{y}\bigg)^n.\label{WWWWWW-999-kkkk}
\end{align}
\end{tl}
\pf Let  $r=1$ and $a_1=b_1$ in \eqref{WWWWWW-9}. In this case, \eqref{WWWWWW-9} reduces to
\begin{align*}
&\sum_{n=0}^\infty \frac{\tau(n)^2}{(q,q/a_2; q)_n}
~{ }_{1} \phi_0\left[\begin{array}{c}
a_{2}q^{-n} \\
-
\end{array}; q, y\right]\bigg(\frac{xq}{ya_2}\bigg)^n\nonumber\\
&=(x;q)_\infty \sum_{n=0}^\infty \frac{\tau(n)^2}{(q,x,q/a_2; q)_n}{ }_{1} \phi_0\left[\begin{array}{c}
a_{2}q^{-n} \\
-
\end{array}; q, yq^n\right]
\bigg(\frac{xq}{ya_2}\bigg)^n.
\end{align*}
Obviously, each ${}_1\phi_0$ series on both sides can be evaluated in closed form by virtue of the $q$-binomial theorem (cf. \cite[(II.3)]{10}). In the sequel, we obtain 
\begin{align*}
\sum_{n=0}^\infty\frac{(q/ya_2;q)_n}{(q,q/a_2;q)_n}\tau(n)
x^n=(x;q)_\infty \sum_{n=0}^\infty \frac{\tau(n)^2(y; q)_n}{(q,x,q/a_2; q)_n}
\bigg(\frac{xq}{ya_2}\bigg)^n.
\end{align*}
Making the replacement $a_2\to q/z$, we obtain \eqref{WWWWWW-999-kkkk}.
\qed
\begin{remark} It is worth pointing out that \eqref{WWWWWW-999-kkkk} turns out to be  \cite[Eq.(6.1.13)]{12-4}.
Observe that the series on the right-hand side of \eqref{WWWWWW-999-kkkk} is symmetric with respect to $x$ and $z$. This leads us to
\begin{align}
\frac{1}{(x;q)_\infty}\sum_{n=0}^\infty\frac{(z/y;q)_n}{(q,z;q)_n}\tau(n)
x^n=\frac{1}{(z;q)_\infty}\sum_{n=0}^\infty\frac{(x/y;q)_n}{(q,x;q)_n}\tau(n)
z^n.
\end{align}
We refer the reader to  \cite[Eq.(6.1.3)]{12-4} for applications of these two transformations in proving $q$-series identities of Rogers-Ramanujan-Slater type.
\end{remark}
In their recent paper \cite{wei}, C. Wei et al's established the following $q$-series expansion via the  analytic method.  Here, by virtue of \eqref{WWWWWW-999-kkkk}, we are able to give a short proof for this expansion.
\begin{lz}[cf. \mbox{\cite[Thm. 1.1]{wei}}]
 \begin{align}(z ; q)_{\infty} \sum_{n=0}^{\infty} \frac{(-x /z ; q)_n}{(q ; q)_n(z ; q)_n} q^{\binom{n}{2}} z^n=\sum_{j, k \geq 0} \frac{q^{j^2+2 j k+2 k^2-j-k}}{(q ; q)_j\left(q^2 ; q^2\right)_k} x^j z^{2 k}.\label{WWWWWW-999-kkkk-1}
\end{align}
\end{lz}
\pf First, replace  $(x,y)$ with $(-z,-z^2/x)$ in \eqref{WWWWWW-999-kkkk} and then multi-sum ply both sides with $(z;q)_\infty$. We obtain
\begin{align}
(z;q)_\infty
\sum_{n=0}^\infty \frac{(-x/z;q)_n}{(q,z;q)_n}q^{\binom{n}{2}}
z^n&=(z^2;q^2)_\infty \sum_{n=0}^\infty \frac{\tau(n)^2(-z^2/x; q)_n}{(q;q)_n(z^2; q^2)_n}x^n\label{surprised}\\
&=\sum_{n\geq j\geq 0}(z^2q^{2n}; q^2)_\infty \frac{q^{n(n-1)+\binom{j}{2}}}{(q;q)_j(q;q)_{n-j}}x^{n-j}z^{2j}\nonumber\\
&=\sum_{n\geq  j\geq 0, k\geq 0} (-1)^{k} \frac{q^{n(n-1)+\binom{j}{2}+k(k-1)}}{(q;q)_j(q;q)_{n-j}(q^2;q^2)_k}x^{n-j}z^{2(j+k)}q^{2nk}.\nonumber
\end{align}
Next, shift $j$ to $n-J.$ So it remains to show
\begin{align*}\sum_{n\geq  J\geq 0, k\geq 0} (-1)^k \frac{q^{n(n-1)+\binom{n-J}{2}+k(k-1)+2nk}}{(q;q)_{n-J}(q;q)_{J}(q^2;q^2)_k}x^{J}z^{2(n-J+k)}=\sum_{j, k \geq 0} \frac{q^{j^2+2 j k+2 k^2-j-k}}{(q ; q)_j\left(q^2 ; q^2\right)_k} x^j z^{2 k}.
\end{align*}
To this end,  one way is to equate the coefficients of the $x^J$ on both sides.  The result is as follows.
\begin{align*}\sum_{n\geq  J\geq 0, k\geq 0} (-1)^k \frac{q^{n(n-1)+\binom{n-J}{2}+k(k-1)+2nk}}{(q;q)_{n-J}(q^2;q^2)_k}z^{2(n-J+k)}
&=\sum_{k=0}^\infty \frac{q^{J^2+2 J k+2 k^2-J-k}}{\left(q^2 ; q^2\right)_k} z^{2 k},
\end{align*}
which is  equivalent to
\begin{align*}
\sum_{k=0}^\infty (-1)^k \frac{q^{k(k-1)+J(J-1)+2kJ}}{(q^2;q^2)_k}z^{2k}
\sum_{N:=n-J=0}^\infty\frac{q^{3\binom{N}{2}+2JN+2kN}}{(q;q)_{N}}z^{2N}&=\sum_{k=0}^\infty \frac{q^{J^2+2 J k+2 k^2-J-k}}{\left(q^2 ; q^2\right)_k} z^{2 k}.
\end{align*}
Clearly, on dividing both sides by $q^{J(J-1)}$,  we come up with
\begin{align*}
\sum_{k=0}^\infty (-1)^k \frac{q^{k(k-1)+2kJ}}{(q^2;q^2)_k}z^{2k}
\sum_{N=0}^\infty\frac{q^{3\binom{N}{2}+2kN}}{(q;q)_{N}}(q^Jz)^{2N}&=\sum_{k=0}^\infty \frac{q^{2 J k+2 k^2-k}}{\left(q^2 ; q^2\right)_k} z^{2 k}.
\end{align*}
By exchanging the  order of sums, we get
\begin{align*}
\sum_{N=0}^\infty\frac{q^{3\binom{N}{2}}}{(q;q)_{N}}(q^{J}z)^{2N}\sum_{k=0}^\infty (-1)^k \frac{q^{k(k-1)}}{(q^2;q^2)_k}(q^{2J+2N}z^2\big)^k&=\sum_{k=0}^\infty \frac{q^{2 J k+2 k^2-k}}{\left(q^2 ; q^2\right)_k} z^{2 k}.
\end{align*}
At this stage, the inner sum on the left side can be evaluated in closed-form by the $q$-binomial theorem. We finally arrive at
\begin{align}
(W^2;q^2\big)_\infty\sum_{N=0}^\infty
\frac{q^{3\binom{N}{2}}}{(q;q)_{N}(W^2;q^2\big)_N}W^{2N}&=\sum_{k=0}^\infty \frac{q^{2 k^2-k}}{\left(q^2 ; q^2\right)_k} W^{2 k}.\label{happy}
\end{align}
In the above, we write $W$ for $q^{J}z$. The validity of \eqref{happy} can be justified
 by the limit $x\to 0$ of \eqref{surprised}. Precisely,  when $x\to 0$,
  it is easy to check
\begin{align*}(z^2;q^2)_\infty \sum_{n=0}^\infty \frac{q^{3\binom{n}{2}}z^{2n}}{(q;q)_n(z^2; q^2)_n}
&=(z;q)_\infty
\sum_{n=0}^\infty \frac{q^{\binom{n}{2}}z^n}{(q,z;q)_n}\\
&=\sum_{n=0}^\infty \frac{\tau(n)z^n}{(q;q)_n}\sum_{k=0}^n\left[n \atop k\right]_q (-1)^k=\sum_{n=0}^\infty \frac{(q;q^2)_{n}}{(q;q)_{2n}}\tau(2n)z^{2n}.\end{align*}
Note that  the last equality is built on the Gauss identity (cf. \cite{gauss}).
 This completes the proof of Wei et al.'s expansion.
\qed

\subsubsection{Transformations with ${}_{2}\phi_1$ series involved}
Continuing in the above way, we may yet deduce from  \eqref{WWWWWW-9}
\begin{tl}For $b_1\neq 1, b_1\neq a_i ~(i=1,2)$, it holds
\begin{align}
&\sum_{n=0}^\infty \frac{q^{n^2}(a_1q/b_1;q)_n}{(q,q/a_1,q/a_2;q)_n}
\bigg(\frac{x}{a_1}\bigg)^n~{ }_{2} \phi_1\left[\begin{array}{c}
q^{-n},a_{2}q^{-n} \\
b_1q^{-n}/a_1
\end{array}; q, \frac{b_1}{a_2}\right]\nonumber\\
&=(x;q)_\infty \sum_{n=0}^\infty \frac{\tau(n)^3}{(q,x ; q)_n}\frac{(b_1/(a_1a_2);q)_n}{(q/a_1,q/a_2;q)_n}
\bigg(\frac{xq^2}{b_1}\bigg)^n.\label{WWWWWW-999}
\end{align}
\end{tl}
\pf For the case $r=1$, \eqref{WWWWWW-9}  is equal to
\begin{align}
&\sum_{n=0}^\infty \frac{\tau(n)^2}{(q ; q)_n}\frac{(q/b_1;q)_n}{(q/a_1,q/a_2;q)_n}
~{ }_{2} \phi_1\left[\begin{array}{c}
a_1q^{-n},a_{2}q^{-n} \\
b_1q^{-n}
\end{array}; q, y\right]
\bigg(\frac{xb_1q}{ya_1a_2}\bigg)^n\label{WWWWWW-999-kk}\\
&=(x;q)_\infty \sum_{n=0}^\infty \frac{\tau(n)^2}{(q,x; q)_n}\frac{(q/b_1;q)_n}{(q/a_1,q/a_2;q)_n}{ }_{2} \phi_1\left[\begin{array}{c}
a_1q^{-n}, a_{2}q^{-n} \\
b_1q^{-n}
\end{array}; q, yq^n\right]
\bigg(\frac{xb_1q}{ya_1a_2}\bigg)^n.\nonumber
\end{align}
Moreover, when $y=b_1/(a_1a_2)$, we may apply the Gauss summation formula \cite[(II.8)]{10} or \eqref{gauss-id} to  find
\[{ }_{2} \phi_1\left[\begin{array}{c}
a_1q^{-n}, a_{2}q^{-n} \\
b_1q^{-n}
\end{array}; q, yq^n\right]\bigg|_{y=b_1/(a_1a_2)}=\frac{(b_1/a_1,b_1/a_2;q)_\infty}{(b_1q^{-n},b_1q^n/(a_1a_2);q)_\infty}.\]
In addition, by Heine's second transformation  of ${}_2\phi_1$ series \cite[(III.2)]{10}, it holds
\[{ }_{2} \phi_1\left[\begin{array}{c}
a_1q^{-n},a_{2}q^{-n} \\
b_1q^{-n}
\end{array}; q, \frac{b_1}{a_1a_2}\right]=\frac{(b_1/a_2,b_1q^{-n}/a_1;q)_\infty}{(b_1q^{-n},b_1/(a_1a_2);q)_\infty}{ }_{2} \phi_1\left[\begin{array}{c}
q^{-n},a_{2}q^{-n} \\
b_1q^{-n}/a_1
\end{array}; q, \frac{b_1}{a_2}\right].\]
Substituting these results into  \eqref{WWWWWW-999-kk} and simplifying further, we  get the assertion of \eqref{WWWWWW-999}.
\qed

From \eqref{WWWWWW-999} we further  deduce
\begin{lz}The following $q$-series transformations hold:
\begin{description}
\item[(i) For $b_1=q$ in \eqref{WWWWWW-999}:]
\begin{align}&\sum_{n=0}^\infty \frac{q^{n^2}(a_1;q)_n}{(q,q/a_1,q/a_2;q)_n}
\bigg(\frac{x}{a_1}\bigg)^n~{ }_{2} \phi_1\left[\begin{array}{c}
q^{-n},a_{2}q^{-n} \\
q^{1-n}/a_1
\end{array}; q, \frac{q}{a_2}\right]\nonumber\\
&=(x;q)_\infty \sum_{n=0}^\infty \frac{\tau(n)^3}{(q,x ; q)_n}\frac{(q/(a_1a_2);q)_n}{(q/a_1,q/a_2;q)_n}
(xq)^n.\label{WWWWWW-99999}
\end{align}
\item[(ii) For $a_1=a=1/a_2$ in \eqref{WWWWWW-99999}:]
\begin{align}&\sum_{n=0}^\infty \frac{q^{n^2}(a;q)_n}{(q,q/a,aq;q)_n}
\bigg(\frac{x}{a}\bigg)^n~{ }_{2} \phi_1\left[\begin{array}{c}
q^{-n},q^{-n}/a \\
q^{1-n}/a
\end{array}; q, aq\right]\nonumber\\ &=(x;q)_\infty \sum_{n=0}^\infty \frac{q^{n(3 n-1)/2}(-x)^n}{(x,aq,q/a;q)_n}.\label{WWWWWW-999999-0}
\end{align}
\item[(iii) For $a=1$ in \eqref{WWWWWW-999999-0}:]
\begin{align}\sum_{n=0}^\infty \frac{q^{n( n-1)/2}(-x)^n}{(q; q)_n^2}=(x;q)_\infty\sum_{n=0}^\infty \frac{q^{n(3 n-1)/2}(-x)^n}{(x,q,q; q)_n}.\label{WWWWWW-9999999}
\end{align}
\end{description}
\end{lz}

Once applying Transformation \cite[(III.32)]{10} to the ${}_2\phi_1$ series on the left-hand side of \eqref{WWWWWW-999},  we find
\begin{tl}For $b_1\neq a_i~(i=1,2)$, it holds
\begin{align}
&\sum_{n=0}^\infty \frac{\tau(n)x^n}{(q,q/a_1;q)_n}
{ }_{2} \phi_1\left[\begin{array}{c}
q^{-n},a_{1}q/b_1 \\
q/a_2
\end{array}; q, \frac{q^{n+1}}{a_1}\right]\nonumber\\
&=(x;q)_\infty \sum_{n=0}^\infty \frac{\tau(n)^3}{(q,x ; q)_n}\frac{(b_1/(a_1a_2);q)_n}{(q/a_1,q/a_2;q)_n}
\bigg(\frac{xq^2}{b_1}\bigg)^n.\label{WWWWWW-999-1}
\end{align}
\end{tl}

Although \eqref{WWWWWW-999} is restricted to  $b_1\neq a_1$,  the limiting case  $b_1\to a_1$ is  of  interest.
\begin{tl}\begin{align}
&\sum_{n=0}^\infty \frac{q^{n^2}(x/a_1)^n}{(q/a_1,q/a_2;q)_n}
\sum_{k=0}^n \frac{ (a_2q^{-n};q)_k}{ (q;q)_k}\bigg(\frac{a_1}{a_2}\bigg)^k\nonumber\\
&=(x;q)_\infty \sum_{n=0}^\infty \frac{\tau(n)^3}{(q,x,q/a_1; q)_n}\frac{1-1/a_2}{1-q^n/a_2}
\bigg(\frac{xq^2}{a_1}\bigg)^n.\label{WWWWWW-99999-xx}
\end{align}
\end{tl}

As is expected,  Identity \eqref{WWWWWW-99999-xx} indeed cover some  useful $q$-series identities.
\begin{lz}The following $q$-series transformations hold:
\begin{description}
\item[(i)]
The limit of \eqref{WWWWWW-99999-xx} as $a_2\to\infty$ yields
\begin{align}
\sum_{n=0}^\infty \frac{q^{n^2}(x/a_1)^n}{(q/a_1;q)_n}
\sum_{k=0}^n \frac{ \tau(k)(a_1q^{-n})^k}{ (q;q)_k}=(x;q)_\infty \sum_{n=0}^\infty \frac{\tau(n)^3}{(q,x,q/a_1; q)_n}
\bigg(\frac{xq^2}{a_1}\bigg)^n.\label{WWWWWW-999-xxx}
\end{align}
\item[(ii)] The limit as $a_1\to 0$ of \eqref{WWWWWW-999-xxx}, in turn, reduces to the special case of \cite[(6.1.4)]{12-4}
\begin{align}
\sum_{n=0}^\infty \tau(n)x^n
=(x;q)_\infty \sum_{n=0}^\infty \frac{q^{n^2}x^n}{(q,x; q)_n}.\label{WWWWWW-999-xxx-00}
\end{align}
\item[(iii)] Alternately, let $a_1=a_2q=q^2/t$ in \eqref{WWWWWW-99999-xx} . Then
\begin{align}
\sum_{n=0}^\infty \frac{q^{n^2-2n}(xt)^n}{(t/q,t;q)_n}
\sum_{k=0}^n \frac{ (q^{1-n}/t;q)_k}{ (q;q)_k}q^k
=(x;q)_\infty\sum_{n=0}^\infty \frac{\tau(n)^3}{(q,x,t; q)_n}
(xt)^n.\label{WWWWWW-999-x}
\end{align}
\item[(iv)] A direct comparison  with the special case that $z=t, y\to\infty$  of \eqref{WWWWWW-999-kkkk} yields 
\begin{align}
\sum_{k=0}^n \frac{ (q^{1-n}/t;q)_k}{ (q;q)_k}q^k
=\frac{(t/q;q)_n}{(q;q)_n}\frac{q^{-n^2/2+3n/2}}{(-t)^n}.\label{WWWWWW-999-kkkk-added}
\end{align}
\end{description}
\end{lz}

\subsubsection{Transformations with ${}_{3}\phi_2$ series involved}
Now we deal with the  $r=2$ case of Theorem \ref{firsttrans}, which takes the following form.
\begin{dl}\label{thirsttrans}
\begin{align}
&\sum_{n=0}^\infty \frac{q^{n^2}}{(q ; q)_n}\frac{(q/b_1,q/b_2;q)_n}{(q/a_1,q/a_2,q/a_{3};q)_n}
\bigg(\frac{xb_1b_2}{ya_1a_2a_{3}}\bigg)^n\nonumber\\
&\qquad\qquad\times{ }_{3} \phi_2\left[\begin{array}{c}
a_1q^{-n}, a_2q^{-n}, a_{3}q^{-n} \\
b_1q^{-n}, b_2q^{-n}
\end{array}; q, y\right]\nonumber\\
&=(x ; q)_{\infty} \sum_{n=0}^\infty \frac{q^{n^2}}{(q,x ; q)_n}\frac{(q/b_1,q/b_2;q)_n}{(q/a_1,q/a_2,q/a_{3};q)_n}
\bigg(\frac{xb_1b_2}{ya_1a_2 a_{3}}\bigg)^n\nonumber\\
&\qquad\qquad\qquad\times{ }_{3} \phi_2\left[\begin{array}{c}
a_1q^{-n}, a_2q^{-n}, a_{3}q^{-n} \\
b_1q^{-n},  b_2q^{-n}
\end{array}; q, yq^n\right].\label{WWWWWW-9-added}
\end{align}
\end{dl}
Recall that
\begin{yl}[\mbox{cf. \cite[(III.9)/(III.10)]{10}}]
\begin{align}{}_3 \phi_2\left[\begin{array}{c}
a, b, c \\
d, e
\end{array}; q, \frac{d e}{a b c}\right]
& =\frac{(e / a, d e / b c ; q)_{\infty}}{(e, d e / a b c ; q)_{\infty}} {}_3 \phi_2\left[\begin{array}{c}
a, d / b, d / c \\
d, d e / b c
\end{array} ; q, \frac{e}{a}\right]\label{phi32-one} \\
& =\frac{(b, d e / a b, d e / b c ; q)_{\infty}}{(d, e, d e / a b c ; q)_{\infty}}{}_3 \phi_2\left[\begin{array}{c}
d / b, e / b, d e / a b c \\
d e / a b, d e / b c
\end{array};q, b\right].\label{phi32-two}
\end{align}
\end{yl}
Upon taking $y=b_1b_2/a_1a_2a_3$ and applying both  \eqref{phi32-one} and \eqref{phi32-two} to the ${}_3\phi_2$ series on the right-hand side of \eqref{WWWWWW-9-added} successively, we get
\begin{tl}\label{thirsttrans-one}
\begin{align}
&\sum_{n=0}^\infty \frac{q^{n^2}x^n}{(q ; q)_n}\frac{(q/b_1,q/b_2;q)_n}{(q/a_1,q/a_2,q/a_{3};q)_n}\nonumber\\
&\qquad\qquad\times{ }_{3} \phi_2\left[\begin{array}{c}
a_1q^{-n}, a_2q^{-n}, a_{3}q^{-n} \\
b_1q^{-n}, b_2q^{-n}
\end{array}; q, \frac{b_1b_2}{a_1a_2a_3}\right]\nonumber\\
&=
\frac{(x,b_2/a_1,b_1b_2/a_2a_3;q)_\infty}{(b_2,b_1 b_2/a_1a_2a_3;q)_\infty} \sum_{n=0}^\infty \frac{\tau(n)^{3}}{(q,x ; q)_n}\frac{(q/b_1,b_1 b_2/a_1a_2a_3;q)_n}{(q/a_1,q/a_2,q/a_{3};q)_n}\bigg(\frac{xq^2}{b_2}\bigg)^n\nonumber\\
&\qquad\qquad\qquad\qquad\qquad\times{ }_{3} \phi_2\left[\begin{array}{c}a_1q^{-n}
,b_1/a_2,b_1/a_3\\
b_1q^{-n},b_1b_2/a_2a_3
\end{array}; q, \frac{b_2}{a_1}\right]\label{WWWWWW-9-added1}\\
&=
\frac{(x,a_2,b_1b_2/a_1a_2,b_1 b_2/a_2a_3;q)_\infty}{(b_1,b_2,b_1 b_2/a_1a_2a_3;q)_\infty} \sum_{n=0}^\infty \frac{\tau(n)^{3}}{(q,x ; q)_n}\frac{(b_1 b_2/a_1a_2a_3;q)_n}{(q/a_1,q/a_{3};q)_n}\bigg(\frac{xa_2q^2}{b_1b_2}\bigg)^n\nonumber\\
&\qquad\qquad\qquad\qquad\qquad\times{ }_{3} \phi_2\left[\begin{array}{c}
b_1/a_2,b_2/a_2,b_1b_2q^{n}/a_1a_2a_3 \\
b_1b_2/a_1a_2,b_1b_2/a_2a_3
\end{array}; q, a_2q^{-n}\right].\label{WWWWWW-9-added2}
\end{align}
\end{tl}
\pf Since, by \eqref{phi32-one},  it holds
\begin{align}
&{ }_{3} \phi_2\left[\begin{array}{c}
a_1q^{-n}, a_2q^{-n}, a_{3}q^{-n} \\
b_1q^{-n},  b_2q^{-n}
\end{array}; q,\frac{b_1b_2q^n}{a_1a_2a_3}\right]\label{WWWWWW-9-added222}\\
&=\frac{(b_2/a_1,b_1b_2/a_2a_3;q)_\infty}{(b_2q^{-n},b_1 b_2q^n/a_1a_2a_3;q)_\infty} { }_{3} \phi_2\left[\begin{array}{c}
a_1q^{-n}, b_1/a_2,b_1/a_{3}\\
b_1q^{-n},b_1  b_2/a_2a_3
\end{array};q, \frac{b_2}{a_1}\right],\nonumber\end{align}
while, from \eqref{phi32-two}, we find
\begin{align}
&{ }_{3} \phi_2\left[\begin{array}{c}
a_1q^{-n}, a_2q^{-n}, a_{3}q^{-n} \\
b_1q^{-n},  b_2q^{-n}
\end{array}; q,\frac{b_1b_2q^n}{a_1a_2a_3}\right]\label{WWWWWW-9-added333}\\
&=
\frac{(a_2q^{-n},b_1b_2/a_1a_2,b_1b_2/a_2a_3;q)_\infty}{(b_1q^{-n},b_2q^{-n},b_1 b_2q^n/a_1a_2a_3;q)_\infty}{ }_{3} \phi_2\left[\begin{array}{c}b_1/a_2,b_2/a_2,b_1b_2q^n/a_1a_2a_3\\
b_1b_2/a_1a_2,b_1b_2/a_2a_3
\end{array}; q, a_2q^{-n}\right].\nonumber\end{align}
Therefore, replacing the ${}_3\phi_2$ series on the right-hand side of \eqref{WWWWWW-9-added} with \eqref{WWWWWW-9-added222} and \eqref{WWWWWW-9-added333} respectively, we get \eqref{WWWWWW-9-added1} and \eqref{WWWWWW-9-added2}.
\qed

The case $b_1=a_2$ of \eqref{WWWWWW-9-added2} is as follows.
\begin{lz}

\begin{align}
&\sum_{n=0}^\infty \frac{q^{n^2}x^n}{(q ; q)_n}\frac{(q/b_2;q)_n}{(q/a_1,q/a_{3};q)_n}{ }_{2} \phi_1\left[\begin{array}{c}
a_1q^{-n},  a_{3}q^{-n} \\
 b_2q^{-n}
\end{array}; q, \frac{b_2}{a_1a_3}\right]\nonumber\\
&=
\frac{(x,b_2/a_1,b_2/a_3;q)_\infty}{(b_2, b_2/a_1a_3;q)_\infty} \sum_{n=0}^\infty \frac{\tau(n)^{3}}{(q,x ; q)_n}\frac{( b_2/a_1a_3;q)_n}{(q/a_1,q/a_{3};q)_n}\bigg(\frac{xq^2}{b_2}\bigg)^n\label{WWWWWW-9-added1-11}.
\end{align}
It is to be noted that   the special cases $a_1=1$ of both  \eqref{WWWWWW-9-added1-11} and \eqref{WWWWWW-999} are equivalent to each other after parameter relabeling.
\end{lz}
We end our discussions with a special instance of \eqref{ssssss-added}.
\begin{lz}
\begin{align}
 &\sum_{n=0}^\infty\frac{\tau(n)(
q/b_1,q/b_2;q)_n}{(q,q,q/a_2,q/a_{3};q)_n}\bigg(\frac{xb_1b_2}{a_2a_3}\bigg)^n\label{ssssss-addeded}\\
&=(x;q)_\infty\sum_{i\geq k\geq  0} \frac{x^iq^{ik}\tau(k)}{(q;q)_i(q,x; q)_k}
\frac{(q/b_1,b_1b_2/a_2a_3;q)_{i-k}}{(q,q/a_2,q/a_{3};q)_{i-k}} { }_3 \phi_2\left[\begin{array}{c}
q^{-(i-k)}, b_1/a_2, b_1/a_3 \\
b_1q^{k-i},b_1b_2/a_2a_3
\end{array} ; q, q\right].\nonumber
\end{align}
\end{lz}
\pf Let $a_1=1, r=2,$ and $xy=q$. Observe  that in this case,  $(a_1;q)_{k-i}\neq 0$ only if $k-i\leq 0$. Thus \eqref{ssssss-added}  reduces to
\begin{align}
 &\sum_{n=0}^\infty\frac{\tau(n)(q/b_1,q/b_2;q)_n}
 {(q,q,q/a_2,q/a_{3};q)_n}\bigg(\frac{xb_1b_2}{a_2a_3}
 \bigg)^n\label{cccccc}\\
&=(x;q)_\infty\sum_{i\geq k\geq  0} \frac{x^iq^{ik+k-i}\tau(k)}{(q;q)_i(q,x; q)_k}
\frac{(1,a_2,a_{3};q)_{k-i}}{(b_1,b_2;q)_{k-i}}{ }_{3} \phi_{2}\left[\begin{array}{c}
q^{-(i-k)},q^{k-i}a_2,q^{k-i}a_{3}\\
q^{k-i}b_1,q^{k-i}b_2
\end{array}; q,q\right].\nonumber
\end{align}
Observe that by \cite[(III.11)]{10}, it holds
\begin{align*}&{ }_{3} \phi_{2}\left[\begin{array}{c}
q^{-(i-k)},q^{k-i}a_2,q^{k-i}a_{3}\\
q^{k-i}b_1,q^{k-i}b_2
\end{array}; q,q\right]\\
&=\frac{(b_1b_2/a_2a_3 ; q)_{i-k}}{(q^{k-i}b_2; q)_{i-k}}\left(\frac{a_2a_3q^{k-i}}{b_1}\right)^{i-k}{ }_3 \phi_2\left[\begin{array}{c}
q^{-(i-k)}, b_1/a_2, b_1/a_3 \\
b_1q^{k-i},b_1b_2/a_2a_3
\end{array} ; q, q\right].
\end{align*}
Again from \eqref{basic33} it follows
\begin{align*}
 \frac{(1,a_2,a_{3};q)_{k-i}}{(b_1,b_2;q)_{k-i}}=
\frac{(q/b_1,q/b_2;q)_{i-k}}{(q,q/a_2,q/a_{3};q)_{i-k}} \bigg(\frac{b_1b_2q}{a_2a_3}\bigg)^{i-k}\tau(i-k).
\end{align*}
Putting  these results together, we may  reduce \eqref{cccccc} to \eqref{ssssss-addeded}. 
\qed


\begin{thebibliography}{99}
\bibliographystyle{amsplain}
\bibitem{9-0}G. E. Andrews, multi-sum ple $q$-series identities, Houston J. Math. 7(1) (1981), 110-122.
\bibitem{9}G. E. Andrews and B. Berndt, Rogers-Ramanujan-Slater type identities. In: Ramanujan's Lost Notebook. Springer, New York, NY. 2005. \url{https://doi.org/10.1007/0-387-28124-X_12}
\bibitem{12-jmaa}Z. Cao and L. Wang, multi-sum -sum   Rogers-Ramanujan type identities, J. Math. Anal. Appl. 522 (2023), Art. 126960.
\bibitem{12} Z. Cao, H. Rosengren, and L. Wang,    On some double Nahm sums of Zagier,  J. Combin. Theory, Series A, Vol. 202 (2024),105819.
\bibitem{10}G. Gasper and M. Rahman, Basic Hypergeometric Series (2nd Edition), Cambridge University Press, Cambridge, 2004.
\bibitem{gauss}C. F. Gauss, Werke, Vol. 2, K$\ddot{o}$niglichen Gesellschaft der Wissenschaften, G$\ddot{o}$ttingen, 1876.
 \bibitem{12-4} J. Mc. Laughlin, A. V. Sills, and P. Zimmer, Rogers-Ramanujan-Slater type identities, Electron J. Comb. 29 (2022), \#P00
\bibitem{name}W. Nahm, Conformal field theory and torsion elements of the Bloch group, in Frontiers in Number Theory, Physics, and Geometry, II, Springer, 2007, 67-132.
 \bibitem{12-3} A.V. Sills,  Identities of the Rogers-Ramanujan-Slater type, Int. J. Number Theory 3(2) (2007), 293-323.
\bibitem{11} M. Vlasenko and S. Zwegers, Nahm's conjecture: asymptotic computations and counterexamples, Commun. Number Theory Phys. 5(3) (2011), 617-642.
\bibitem{12-2} L. Wang, Identities on Zagier's rank two examples for Nahm's conjecture, arXiv:2210.10748.
\bibitem{wei}C. Wei, Y. Yu, and G. Ruan, Multi-sum dimensional Rogers-Ramanujan type identities with parameters, 	arXiv:2302.00357.
   \bibitem{13} D. Zagier, The dilogarithm function, in Frontiers in Number Theory, Physics and Geometry, II, Springer, 2007, 3-65.
\end{thebibliography}
\end{document}